\documentclass{amsart}
\theoremstyle{plain}
\newtheorem{thm}{Theorem}[section]

\newtheorem{cor}[thm]{Corollary}
\newtheorem{conj}[thm]{Conjecture}
\newtheorem{prob}[thm]{Problem}

\theoremstyle{definition}
\newtheorem{defn}[thm]{Definition}

\begin{document}

\title{The Bing-Borsuk and the Busemann Conjectures}
\author{Denise M. Halverson and Du\v san Repov\v s}

\address{Department of Mathematics, Brigham Young University,
Provo, UT 84602}
\email{deniseh@math.byu.edu}
\address{Faculty of Mathematics and Physics, and 
Faculty of Education,
University of Ljubljana, 
P. O. B. 2964,
Ljubljana, Slovenia 1001}
\email{dusan.repovs@guest.arnes.si}



\begin{abstract}
We present two classical conjectures concerning the characterization
of manifolds: the Bing Borsuk Conjecture asserts that every
$n$-dimensional homogeneous ANR is a topological $n$-manifold,
whereas the Busemann Conjecture asserts that every $n$-dimensional
$G$-space is a topological $n$-manifold. The key object in both
cases are so-called {\it generalized manifolds}, i.e. ENR homology
manifolds. We look at the history, from the early beginnings
to the present  day. We also list several
open problems and related conjectures.
\end{abstract}

\date{\today}

\maketitle

\section{Introduction}

In this paper we will survey two famous conjectures, their
relationship to each other, as well as their relationship to other
famous unsolved problems.  The first, and most general, is the
Bing-Borsuk Conjecture \cite{Bing-Borsuk} which states that any
$n$-dimensional homogeneous ENR space is an $n$-manifold.  
The $3$-dimensional
Bing-Borsuk Conjecture 
implies the celebrated Poincar\'{e}
Conjecture, recently proven by Perel'man \cite{Morgan-Tian}.
Given the complexity of the proof the Poincar\'{e} Conjecture, it is
understandable why in particular, 
the $3$-dimensional Bing-Borsuk Conjecture remains unsolved.

The Busemann Conjecture \cite{Busemann1}-\cite{Busemann3} is a special case of the
Bing-Borsuk Conjecture and it states that Busemann $G$-spaces are
manifolds.  Busemann $G$-spaces are well known to be homogeneous.
Therefore the truth of the Bing-Borsuk Conjecture immediately
implies the truth of the Busemann Conjecture.  The Busemann
Conjecture has been proven for $G$-spaces of dimensions $\leq 4$
(see \cite{Busemann2}\cite{Krakus}\cite{Thurston}).

Three other famous conjectures related to the
Busemann Conjecture are the
Moore conjecture  and
the two de Groot conjectures.
The Moore conjecture
\cite{DavermanBook}\cite{Wilder}
is to
determine whether or not all resolvable generalized manifolds are
codimension one manifold factors.  Small metric balls in Busemann
$G$-spaces are known to be cones over
their boundaries, and hence have a
local product structure with respect to their boundaries which are
known to be generalized $(n-1)$-manifolds. An affirmative answer to
the Moore conjecture together with the resolvability of Busemann
$G$-spaces implies an affirmative answer to the Busemann Conjecture.

The de Groot Conjectures \cite{de Groot} are to determine whether or
not all absolute cones are balls (resp. whether or not all absolute
suspensions are spheres). Small metric balls in Busemann $G$-spaces
are also absolute cones. It has recently been proven that absolute
cones are $n$-cells in dimensions $n \leq 4$, but there are
counter-examples in higher dimensions \cite{Guilbault}.  The
solution to the $n=4$ case relies upon the Poincar\'{e} Conjecture.
Unfortunately, this result does not provide a solution to the
Busemann Conjecture in dimension $n \geq 5$.

The purpose of this paper is to survey the work that has been done
on these manifold recognition problems.  In Section 2 we will
delineate important properties that are known to be satisfied by
manifolds. In Section 3 we will provide an overview of progress that
has been made towards resolving the Bing-Borsuk Conjecture.  In
Section 4 we will do the same for the Busemann Conjecture.  In
Section 5 we will discuss three related problems: the Moore Conjecture
and two de Groot Conjectures. In Sections 6 and 7 we will provide a list of
relevant problems that remain unsolved.

\section{Manifolds and Manifold Properties}

An {\it $n$-manifold} is a separable metric space such that each
point has a neighborhood homeomorphic to  the Euclidean $n$-space $\mathbb R^n$.
Although this is a simple definition to state, applying the
definition to verify that a space is a manifold can be a difficult
task.  Thus it is desirable to find alternate methods of detecting
manifolds. In this section we describe the properties and
characteristics known to be possed by manifolds.  The question is
which property or combination of properties are sufficient to imply
that a space is a manifold (see \cite{Cannon}\cite{Repovs1}\cite{Repovs3}).

A topological space $X$ is said to be {\it homogeneous} if
for any
two points $x_{1}, x_{2} \in X$, there is a homeomorphism of $X$
onto itself taking $x_{1}$ to $x_{2}$.  It is a classical result
that closed (i.e. connected compact without boundary) manifolds are homogeneous.
It is the Bing-Borsuk
Conjecture that asks whether a homogeneous space is necessarily a
manifold \cite{Bing-Borsuk}.

A topological space $X$ is said to have the {\it invariance of
domain property} if for every pair of homeomorphic subsets
$U,V\subset X$, $U$ is open if and only if $V$ is open. Brouwer
\cite{Brouwer1}\cite{Brouwer2} proved a century ago
that every topological $n$-manifold
has the invariance of domain property. Unfortunately, the invariance
of domain property is not sufficient by itself to characterize
manifolds (see e.g. \cite{Duda}\cite{vanDalen}\cite{vanDouwen}\cite{vanMill}).

An  $n$-dimensional compact metric space $X$ is called an {\it
$n$-dimensional Cantor manifold} if whenever $X$ can be expressed as
the union $X = X_1 \cup X_2$  of its proper closed subsets, then
$\text{dim} (X_1 \cap X_2) \geq n-1.$
Urysohn \cite{Urysohn}\cite{Urysohn1}, who introduced this notion in 1925,
proved that every topological $n$-manifold is a Cantor
$n$-manifold. More fundamental results were established by
Aleksandrov \cite{Aleksandrov} in 1928.
Krupski \cite{Krupski} proved in 1993 a more general
result, namely, that every generalized $n$-manifold is a Cantor
$n$-manifold.

A metric space $(X,\rho)$ is said to have the {\it disjoint
$(k,m)$-cells property} ($k,m\in \mathbb{N}$) if for each pair of maps
$f:B^k \to X$ and $g:B^m \to X$ and every $\varepsilon
>0$ there exist maps $f':B^k \to X$ and $g':B^m \to X$ such that
$$\rho(f,f')<\varepsilon, \ \ \hbox{dist}(g,g')<\varepsilon \ \
\hbox{and} \ \ f'(B^k) \cap g'(B^m) = \emptyset.$$ It is well known
that topological manifolds of dimension $n$ have the disjoint
$(k,m)$-cells property for $k+m+1 \leq n$ (see
\cite{Rourke-Sanderson}).  The disjoint $(2,2)$-cells property is
often referred to the {\it disjoint disks property} and plays a key
role in characterizing manifolds of dimension $n \geq 5$.

A space $X$ is said to be {\it locally $k$-connected, $LC^{k}$} ($k
\ge 0$),  if for every point $x\in X$ and every neighborhood
$U\subset X$ of $x$, there exists a neighborhood $V\subset U$ of $x$
such that the inclusion-induced homeomorphisms $\pi_{i\le k}(V) \to
\pi_{i\le k}(U)$ are trivial. Clearly, locally contractible spaces,
such as manifolds and polyhedra, are locally $k$ connected for all
$k$.

Let $Y$ be a metric space. Then $Y$ is said to be an {\it absolute
neighborhood retract (ANR)} provided for every closed embedding $e:Y
\to Z$ of $Y$ into a metric space $Z$, there is an open neighborhood
$U$ of the image $e(Y)$ which retracts to $e(Y)$. That is, there is
a continuous surjection $r:U \to e(Y)$ with $r(x) = x$ for all $x
\in e(Y)$. It is a classical result that finite-dimensional spaces
that are ANR's are characterized as the locally contractible
separable metric spaces \cite{Borsuk}. A {\it Euclidean neighborhood retract}
(ENR) is a finite-dimensional, locally compact, locally contractible
subset $X$ of the Euclidean $n$-space $\mathbb{R}^n$.  It follows
immediately from local contractibility that every topological
manifold is an ENR (hence an ANR).

An $n$-dimensional ($n\in \mathbb{N}$) locally compact Hausdorff
space $X$ is called a {\it $\mathbb{Z}$-homology $n$-manifold}
($n$-hm$_{\mathbb{Z}}$) if for every point $x\in X$ and all $k\in
\mathbb{N}$, $H_k(X,X \setminus \{x\};\mathbb{Z}) \cong
H_k(\mathbb{R}^n,\mathbb{R}^n \setminus \{0\};\mathbb{Z})$.  
Trivially, every
topological manifold is a homology manifold. An $n$-dimensional
topological space $X$ is called a {\it generalized $n$-manifold}
($n\in \mathbb{N}$) if $X$ is an ENR $\mathbb{Z}$-homology
$n$-manifold. It follows that every topological $n$-manifold is a
generalized $n$-manifold. Every generalized $(n\le 2)$-manifold is
known to be a topological $n$-manifold \cite{Wilder}. On the other hand, for every
$n\ge 3$ there exist {\it totally singular} generalized
$n$-manifolds $X$, i.e. $X$ is not locally Euclidean at any point
(see \cite{Cannon}
\cite{DavermanBook}
\cite{Repovs1}
\cite{Repovs3}).

A natural way in which a generalized manifold may arise is as the
image of a cell-like map defined on a manifold. A proper onto map
$f: M \to X$ is said to be {\it cell-like}  
if for every point $x \in
X$, the point-inverse $f^{-1}(x)$ contracts in any neighborhood of
itself (i.e., $f^{-1}(x)$ has the shape of a point)\cite{Lacher}. A space $X$
that is the proper image of a cell-like map is said to be {\it
resolvable}. Trivially every topological manifold is resolvable.

The following classical result attests to the crucial importance of
cell-like maps in geometric topology (it was proved for $n\le 2$ by Wilder
\cite{Wilder}, 
for $n=3$ by Armentrout
\cite{Armentrout},
for $n=4$ by Quinn
\cite{Quinn1},
and for $n\ge 5$ by Siebenmann
\cite{Siebenmann}):

\begin{thm}[Cell-like Approximation Theorem]\label{CAT}
For every $\varepsilon >0$, every $n\in \mathbb N$,
and every cell-like map
$f:M^{n}\to N^{n}$ between topological $n$-manifolds
$M^{n}$ and $N^{n}$,
there exists
a homeomorphism  $h:M^{n}\to N^{n}$ such that
$d(f(x),h(x))<\varepsilon$ for every $x\in M^{n}$.
\end{thm}

The fact that not all resolvable 
generalized
manifolds are manifolds has been
known since the mid 1950's, when Bing \cite{Bing1}\cite{Bing2} constructed his
famous Dogbone space as the cell-like image of a map defined on
$\mathbb R^3$.  Generalized manifolds have been 
a subject of intense studies since 1960's \cite{Repovs3}.
In the mid 1970's Cannon recognized that the disjoint
$(2,2)$-cells property, often referred to as the disjoint disks
property (DDP), plays a key role in characterizing manifolds of
dimension $n \geq 5$.
Recall that a metric space $X$ is said to have the {\it disjoint disks property}
(DDP) if for every $\varepsilon >0$ and every pair of maps $f,g:B^2
\to X$ there exist $\varepsilon$-approximations $f', g':B^2 \to X$
with disjoint images $f'(B^2)\cap g'(B^2)=\emptyset$.

Cannon
\cite{Cannon2}
utilized the DDP property
to solve the celebrated Double suspension problem
\cite{Cannon}
\cite{Latour}
\cite{Repovs1}
which asks if the double suspension
$\Sigma^2(H^n)$ of an arbitrary homology $n$-sphere
$H^n$, $n\ge 3$, is the $(n+2)$-sphere $S^{n+2}$.
In 1978, Edwards
\cite{E}
\cite{Edwards}
generalized Cannon's results by proving the famous
Cell-like approximation theorem (for a detailed proof for $n=5$ see
\cite{Daverman-Halverson} and for $n\ge 6$ see
\cite{DavermanBook}),
and at the same time also gave an affirmative answer of
the high-dimensional case of the
Manifold Recognition Problem (which asks if every
resolvable generalized $(n\ge 3)$-manifold with with the
"appropriate amount of general position"
is a topological $n$-manifold
\cite{Cannon}
\cite{Cannon1}
\cite{CRT}
\cite{Daverman1}
\cite{Repovs1}-\cite{Repovs3}):

\begin{thm}[Edwards]\label{E}
For $n \geq 5$, topological
$n$-manifolds are precisely the $n$-dimensional resolvable spaces
with the disjoint disks property.
\end{thm}

An analogous result for $3$-manifolds was proved in the early 1980's
by Daverman and Repov\v s \cite{DaR1} \cite{DaR2} (whereas only
partial results are known in dimension $4$, see \cite{BDVW}
\cite{DaR1}). A metric space $X$ is said to have the {\it Spherical
simplicial approximation property} (SSAP) if for each $\mu:S^{2}\to
X$ and each $\varepsilon >0,$ there exist a map $\psi:S^{2}\to X$
and a finite topological 2--complex $K_{\psi}\subset X$ such that:
(1) $d(\psi,\mu)< \varepsilon$; (2) $\psi(S^{2})\subset K_{\psi}$;
and (3) $X\backslash K_{\psi}$ is 1--FLG in $X$. (The 1--FLG
condition is known to characterize tamely embedded 2--complexes
$K_{\psi}$ in 3--manifolds $M^3$, see \cite{Nicholson}).

\begin{thm}[Daverman-Repov\v s]\label{DR}
Topological
$3$-manifolds are precisely the $3$-dimensional resolvable spaces
with the simplicial spherical approximation property.
\end{thm}

It had  been a long time a question as to whether all generalized
manifolds are resolvable - this was the famous Resolution Conjecture
\cite{BH}-\cite{BP}
\cite{Cannon}
\cite{CBL}-\cite{CHR3}
\cite{HeR}
\cite{Mitchell-Repovs}
\cite{Repovs1}-\cite{Repovs3}:

\begin{conj}[Generalized Manifolds Resolution Conjecture]\label{GMRC}
Every generalized $(n\ge 3)$-manifold has a resolution.
\end{conj}

In dimension 3, the
Generalized Manifolds Resolution Conjecture~\ref{GMRC} implies the
Poincar\'{e} Conjecture \cite{Repovs1} and only special cases are known
\cite{Brin}
\cite{BM}
\cite{BL1}
\cite{BL2}
\cite{CRT}
\cite{DaT}
\cite{Repovs0}
\cite{Repovs}
\cite{RL}
\cite{T1}-\cite{T3}.
In higher dimensions the
Generalized Manifolds Resolution Conjecture~\ref{GMRC}
turns out to be false.  By the results of Bryant et al.
\cite{BFMW} from 1996 which provide the construction
(together the work of Pedersen et al. \cite{PQR} from 2003
which provided some key details on the surgery exact sequence used in the
original construction), it is now known that there exist non-resolvable
generalized $n$-manifolds, for every $n\ge 6$.  In 2007 Bryant et al.
\cite{BFMW2} further strengthened their
result to the following  DDP Theorem:

\begin{thm}[Bryant-Ferry-Mio-Weinberger]\label{BFMW}
There exist non-resolvable generalized $n$-manifolds with the
disjoint disks property, for every $n\ge 7$.
\end{thm}

Hence, generalized manifolds may posses nice general
position properties. Moreover, Krupski has shown that all
generalized manifolds are Cantor manifolds (see Proposition 1.7 of
\cite{Krupski}). Thus, the majority of properties listed above are
known to be insufficient by themselves to characterize manifolds.
Homogeneity is the remaining single candidate.  Is this property
strong enough to characterize manifolds?  Are there other
combinations of these properties that characterize manifolds which
have not yet been discovered? \cite{Repovs3}

\section{The Bing-Borsuk Conjecture}

Bing and Borsuk \cite{Bing-Borsuk} proved in 1965 that for $n<3$, every $n$-dimensional
homogeneous ANR is a topological $n$-manifold.
They also conjectured that this holds in all dimensions:

\begin{conj}[Bing-Borsuk Conjecture]\label{BBC}
Every $n$-dimensional, $n\in \mathbb N$,
homogeneous ANR is a topological $n$-manifold.
\end{conj}

Jakobsche \cite{Jakobsche} proved in 1978 that in dimension $n=3$, the
Bing-Borsuk Conjecture~\ref{BBC}
implies the Poincar\'{e} Conjecture
(see also \cite{Jakobsche-Repovs}).
Given the difficulty of the proof of the
Poincar\'{e} Conjecture \cite{Morgan-Tian}, it is understandable why
the Bing-Borsuk Conjecture~\ref{BBC}
remains unsolved.

\subsection{Partial Results}

Although there is much work to be done before the Bing-Borsuk
Conjecture~\ref{BBC}
will be solved, there are several partial results. In
1970 Bredon \cite{Bredon}\cite{Bredon1} showed the following:

\begin{thm}[Bredon]\label{Bredon}
If $X$ is an $n$-dimensional homogeneous ENR ($n\in \mathbb{N}$) and
for some (and, hence all) points $x\in X$, the groups
$H_k(X,X-\{x\};\mathbb{Z})$ are finitely generated, then $X$ is a
$\mathbb{Z}$-homology $n$-manifold.
\end{thm}

  This theorem was reproved by Bryant \cite{Bryant2} in 1987
with a more geometric argument.
In 1976 Lysko \cite{Lysko} showed:

\begin{thm}[Lysko]\label{L}
Let $X$ be a connected finite-dimensional homogeneous ANR-space.
Then $X$ is a Cantor manifold and possesses the invariance
of domain property.
\end{thm}

In 1985 Seidel \cite{Seidel} proved a similar result in
the case of locally compact, locally homogeneous separable ANR's.

Next, we quote the following result by Krupski \cite{Krupski} from
1993:

\begin{thm}[Krupski]\label{K}
Let $X$ be a homogeneous locally compact space. Then:
(1) If $X$ is an ANR of dimension $> 2$ then $X$ has the
disjoint $(0,2)$-cells property.
(2) If dim $X = n > 0$, $X$ has the disjoint $(0,n-1)$-cells property and
$X$ is an LC$^{n-1}$-space, then local homologies
satisfy $H_k(X, X-\{x\})=0$ for $k<n$ and $H_n(X,X-\{x\}) \ne 0$.
\end{thm}

A topological space $Y$ is said to be {\it acyclic} in dimension 
$n\in \mathbb N$ if $\hbox{\v{H}}^{n}(Y;\mathbb Z)=0$. 
In 2003 Yokoi \cite{Yokoi} established the following algebraic property of $n$-dimensional homogeneous ANR's which is also possessed by topological $n$-manifolds:

\begin{thm}[Yokoi]\label{Y}
Let $X$ be an $n$-dimensional homogeneous ANR continuum which is cyclic in dimension $n$.
Then no compact subset of $X$, acyclic in dimension $n-1$, separates $X$.
\end{thm}

These partial results, demonstrating that homogeneity
implies several of the other manifold properties, indicate why the Bing-Borsuk Conjecture~\ref{BBC} could be true.

\subsection{A special case}

In 1996
Repov\v{s} et al. \cite{RSS} (see  \cite{RSS1} for a 
very geometric proof of the 
$2$-dimensional case) proved the following result
which in some sense can be considered as
a {\it smooth version} of the Bing-Borsuk Conjecture~\ref{BBC}.  Recall that a
subset $K\subset \mathbb{R}^n$ is said to be $C^1$--{\it
homogeneous} if for every pair of points $x, y\in K$ there exist
neighborhoods $O_x, O_y\subset \mathbb{R}^n$ of $x$ and $y$,
respectively, and a $C^1$--{\it diffeomorphism} $$h : (O_x, O_x\cap
K, x)\to (O_y, O_y\cap K, y),$$ i.e. $h$ and $h^{-1}$ have
continuous first derivatives.

\begin{thm}[Repov\v s-Skopenkov-\v S\v cepin]\label{RSS}
Let $K$ be a locally compact (possibly nonclosed) subset of
$\mathbb{R}^n$.
 Then $K$ is $C^1$-homogeneous if and only if $K$ is a
$C^1$--submanifold of $\mathbb{R}^n$.
\end{thm}

This theorem clearly does not work for arbitrary {\it
homeomorphisms}, a counterexample is the {\it Antoine Necklace}
\cite{Antoine} - a wild Cantor set in $\mathbb{R}^3$ which is
clearly {\it homogeneously} (but not $C^1$-{\em homogeneously})
embedded in $\mathbb{R}^3$.
In fact, Theorem~\ref{RSS}
does not even work for {\it
Lipschitz} homeomorphisms, i.e. the maps for which $d(f(x),f(y))<
\lambda \ d(x,y)$, for all $x,y\in X$.
Namely, Male\v si\v c and
Repov\v s \cite{Malesic-Repovs}
proved in 1999 that there exists a Lipschitz homogeneous
wild Cantor set in $\mathbb{R}^3$. Their result was later strengthened by
Garity et al. \cite{GRZ}:

\begin{thm}[Garity-Repov\v s-\v Zeljko]\label{GRZ}
There exist uncountably many rigid Lipschitz homogeneous wild Cantor
sets in $\mathbb{R}^3$.
\end{thm}

\subsection{Alternate Statement}

Daverman and Husch \cite{Daverman-Husch} were able to
determine an equivalent conjecture to the Bing-Borsuk Conjecture. In
order to state this conjecture, recall that a surjective map $p:E
\to B$ between locally compact, separable metric ANR's $E$ and $B$
is said to be an {\it approximate fibration} if $p$ has the {\it
approximate homotopy lifting property} for every space $X$.
Equivalently, whenever $h:X \times I \to B$ and $H:X \times \{0\}
\to E$ are maps such that $p\circ H=h \mid X \times \{0\}$ and
$\varepsilon$ is a cover of $B$, $h$ extends to a map $H:X \times I
\to E$ such that $h$ and $p\circ H$ are $\varepsilon$-close. The
alternate statement of the Bing-Borsuk Conjecture~\ref{BBC} is given as
follows \cite{Bryant1}\cite{Daverman-Husch}\cite{Miller}\cite{West}:

\begin{conj}[Alternate statement of the the Bing-Borsuk Conjecture]\label{ABBC}
Suppose that $X$ is nicely embedded in $\mathbb{R}^{m+n}$, for some
$m \geq 3$, so that $X$ has a mapping cylinder neighborhood
$N=C_\phi$ of a map $\phi:\partial N \to X$, with mapping cylinder
projection $\pi:N \to X$. Then $\pi: N \to X$ is an approximate
fibration.
\end{conj}

\subsection{Modified Bing-Borsuk Conjecture}

Recall that for $n \geq 6$, Bryant et al.
\cite{BFMW} proved in 1996 that there exist non-resolvable
generalized $n$-manifolds, for every $n\ge 6$. Based on earlier work
by Quinn \cite{Quinn2}, these nonresolvable generalized manifolds
must be {\it totally singular}, i.e., have no points with Euclidean
neighborhoods (we may assume these examples are connected).
Moreover, in 2007 Bryant et al.
\cite{BFMW2} strengthened their result to show that there exist
non-resolvable generalized $n$-manifolds with the disjoint disks
property, for every $n\ge 7$.  Based on these results the following
conjecture was proposed:

\begin{conj}[Bryant-Ferry-Mio-Weinberger]\label{BFMW1} Every
generalized $n$-manifold ($n\ge 7$) satisfying  the disjoint disks property, is
homogeneous.
\end{conj}

\indent Note that if Conjecture~\ref{BFMW1} is true, then the
Bing-Borsuk Conjecture~\ref{BBC} is false for $n\ge 7$. In 2002
Bryant \cite{Bryant1} suggested the following modified Bing-Borsuk
Conjecture:

\begin{conj}[Modified Bing-Borsuk Conjecture]\label{MBBC}
Every homogeneous $(n\ge 3)$-dimensional ENR  is a generalized
$n$-manifold.
\end{conj}

A further modification was posed by Quinn \cite{Quinn3} at the 2003
Oberwolfach workshop on exotic homology manifolds.  It is based on
a perturbation of the homogeneity property itself. A space $X$ is
{\it homologically arc-homogeneous} provided that for every path
$\alpha:[0,1] \to X$, the inclusion induced map
$$H_*(X \times 0, X \times 0 - (\alpha(0),0))
\to H_*(X \times I, X \times I - \Gamma(\alpha))$$ is an
isomorphism, where $\Gamma(\alpha)$ denotes the graph of $\alpha$.
The following is the conjecture proposed  by Quinn \cite{Quinn3}
which was proved in 2006 by Bryant \cite{Bryant3}.

\begin{thm}[Bryant]\label{Bryant}
Every $n$-dimensional homologically arc-homogeneous ENR is a generalized
$n$-manifold.
\end{thm}

  This is arguably the strongest result so far relating to
the Bing-Borsuk Conjecture~\ref{BBC}.

\section{The Busemann Conjecture}

The Busemann Conjecture is also a manifold recognition problem, and
is
in fact a special case of the Bing-Borsuk Conjecture~\ref{BBC}. Beginning in
1942, Herbert Busemann \cite{Busemann1}\cite{Busemann2}
developed the notion of a
$G$-space as a way of putting a Riemannian like geometry on a metric
space (and also in an attempt to obtain a "synthetic description" of Finsler's spaces
\cite{Finsler}). A Busemann $G$-space is a metric space that satisfies four
basic axioms on a metric space. These axioms imply the existence of
geodesics, local uniqueness of geodesics, and local extension
properties. These axioms also infer homogeneity and a cone structure
for small metric balls. In 1943, Busemann \cite{Busemann2} proved:

\begin{thm}[Busemann]\label{Busemann}
Busemann $G$-spaces of dimension $n=1,2$ are manifolds.
\end{thm}

  Busemann then proposed the following conjecture \cite{Busemann3}:

\begin{conj}[Busemann Conjecture]\label{BC}
Every $n$-dimensional $G$-space
($n \in \mathbb{N}$) is a topological $n$-manifold.
\end{conj}

  When Busemann \cite{Busemann3} proposed Conjecture~\ref{BC} in
1955, he predicted: {\it Although this (the Busemann Conjecture) is
probably true for any $G$-space, the proof, if the conjecture is
correct, seems quite inaccessible in the present state of topology.}
As we shall see, this prediction proved true.

\subsection{Definitions}

We now formally define Busemann $G$-spaces and state several of the
classical properties of Busemann $G$-spaces.

\begin{defn}
Let $(X,d)$ be a metric space.  $X$ is said to be a {\it Busemann
$G$-space} provided it satisfies the following Axioms of Busemann:

{(i) \emph{Menger Convexity:}} {Given distinct points $x,y \in X$,
there is a point $z \in X-\{x,y\}$ so that $d(x,z) + d(z,y) =
d(x,y).$}

{(ii) \emph{Finite Compactness}:} {Every $d$-bounded infinite set
has an accumulation points.}

{(iii) \emph{Local Extendibility}:} {To every $w \in X$, there is a
positive radius $\rho_w$, such that for any pair of distinct $x,y
\in B(w,\rho_w)$, there is $z \in \text{int } B(w,\rho_w)-\{x,y\}$
such that $d(x,y) + d(y,z) = d(x,z).$}

{(iv) \emph{Uniqueness of the Extension}:} {Given distinct $x,y \in
X$, if there are points $z_1, z_2 \in X$ for which both
$$d(x,y) + d(y,z_i) = d(x,z_i) \quad \text{for } i=1,2, $$
and $$d(y,z_1) = d(y,z_2)$$ hold, then $z_1 = z_2$.}
\end{defn}

From these basic properties, a rich structure on a $G$-space can be
derived. If $(X,d)$ be a $G$-space and $x \in X$, then $(X,d)$
satisfies the following properties:

\begin{itemize}
  \item {\it Complete Inner Metric:} $(X,d)$ is a complete inner metric
space which is locally compact.
  \item {\it Existence of Geodesics:} Any two points in $X$ are joined by a geodesic segment.
  \item {\it Local Uniqueness of Joins:}  There is a positive radius $r_x$ such that any two points $y,z \in B_{r_x}(x)$ in the closed ball are joined by a
unique segment in $X$.
  \item {\it Local Cones:} There is a positive radius $\epsilon_x$ for which the closed metric ball $B_{\epsilon_x}(x)$ is homeomorphic to the cone over its
boundary. That is, $ B_{\epsilon_x} \cong \Delta(S_{\epsilon_x}(x))$
where $S_{\epsilon_x}(x)$ denotes the metric sphere about $x$.
  \item {\it Homogeneity:}  Every $G$-space is homogeneous.  Moreover,
homogeneity homeomorphisms can be chosen so that each is isotopic to
the identity.
\end{itemize}

  It is this last property that makes the Busemann
Conjecture~\ref{BC} a special case of the Bing-Borsuk Conjecture~\ref{BBC}.
The truth
of the Bing-Borsuk Conjecture would imply the truth of the Busemann
Conjecture. Equivalently, if examples of non-manifold Busemann
$G$-spaces could be constructed, the Bing-Borsuk Conjecture would be
settled in the negative.

\subsection{Results in Higher Dimensions}

The first success in resolving the Busemann Conjecture~\ref{BC} in higher
dimensions occured in 1968 when Krakus \cite{Krakus} proved it in dimension $n=3$.

\begin{thm}[Krakus]\label{Krakus}
Buseman $G$-spaces of dimension $n=3$ are manifolds.
\end{thm}

  Krakus applied Borsuk's  $2$-sphere recognition
  criterion
  \cite{Borsuk1}
  to show
that small metric spheres in $3$-dimensional Busemann spaces are
topological $2$-spheres. The truth of the conjecture follows
immediately from the local product structure on small metric spheres
and homogeneity. Unfortunately, this strategy cannot be extended
to higher dimensions due to the lack of similar characterizations of
topological spheres in higher dimensions.

Starting in dimension $n=4$, it can now be seen that Busemann's
prediction of the difficulty of the problem was remarkably accurate.
For example, the case $n=4$ required several modern results and
techniques including sheaf theory \cite{Bredon},
epsilon surgery \cite{CHR3},
resolution
theory \cite{Quinn2},
decomposition theory \cite{DavermanBook},
and theory of $4$-manifolds \cite{Freedman}\cite{Freedman-Quinn}\cite{Scorpan}.
The major breakthrough in dimension $n=4$ and partial results
applicable to higher dimensions were made by Thurston
\cite{Thurston} in 1996:

\begin{thm}[Thurston]\label{Thurston}
Busemann $G$-spaces of dimension
$n=4$ are manifolds. Moreover, every finite dimensional $G$-space is
a generalized $n$-manifold. More precisely, let $(X,d)$ be a
$G$-space, dim$\ X=n<\infty$.  Then for all sufficiently small $r>0$
and $x \in X$, $B_r(x)$ is a homology $n$-manifold with boundary
$\partial B_r(x) = S_r(x)$ and $S_r(x)$ is a homology
$(n-1)$-manifold with empty boundary.
\end{thm}

In 2002 Berestovskii \cite{Berestovskii4} proved the special case of
the Busemann Conjecture~\ref{BC}
for Busemann $G$-spaces that have Alexandrov
curvature bounded above. A Busemann $G$-space $(X,d)$ has Alexandrov
curvature $\leq K$ if geodesic triangles in $X$ are at most as "fat"
as corresponding triangles in a surface $S_K$ of constant curvature
$K$,  i.e., the length of a bisector of the triangle in $X$ is at
most the length of the corresponding bisector of the corresponding
triangle in $S_K$. For example, the boundary of a convex region in
$\mathbb{R}^n$ has a nonnegative Alexandrov curvature (see also
\cite{ABN}
\cite{Berestovskii1}-\cite{Berestovskii3}
\cite{Berestovskii5}
\cite{BGP}
\cite{Hotchkiss}
\cite{Papadopoulos}
\cite{Pogorelov}).

\begin{thm}[Berestovskii]\label{Berestovskii}
Busemann $G$-spaces of dimension $n\geq 5$ having bounded
Aleksandrov curvature bounded above are $n$-manifolds.
\end{thm}

The general case of the Busemann Conjecture~\ref{BC} for $n \geq 5$
remains unsolved: there are many Busemann $G$-spaces which do not
satisfy the condition of Aleksandrov curvature bounded from above
(or below). However, all such examples are known to be topological
manifolds. The simplest example of Busemann G-space which is not
"covered" by Berestovskii's  proof \cite{Berestovskii4} is the
finite-dimensional vector space $(\mathbb{R}^n,| \cdot |)$ in
which the closed unit ball $\{x\in \mathbb{R}^n: |x|\leq 1\}$ is a
strongly convex centrally symmetric convex body in $\mathbb{R}^n$
which is not an ellipse.

\begin{conj}[Busemann $G$-Spaces Resolution Conjecture]\label{BSRC}
Every $(n\ge 5)$-dimensional
Busemann $G$-space has a resolution.
\end{conj}

\section{Related Problems}

Related to the Busemann Conjecture~\ref{BC} are three other famous problems:
the two de Groot Conjectures and the Moore Conjecture.

\subsection{The de Groot Conjectures}

The de Groot Conjectures are two manifold recognition problems for spaces
that are absolute suspensions or absolute cones.  A compact
finite-dimensional metric space $X$ is called an {\it absolute
suspension} (AS) if it is a suspension with respect to any pair of
distinct points  and is called an {\it absolute cone} if it is a
cone with respect to any point. Any space topologically equivalent
to $S^n$ is an absolute suspension and any space topologically
equivalent to $B^n$ is an absolute cone. The question is whether the
converse statements are true. At the 1971 Prague Symposium, de Groot
\cite{de Groot} made the following two conjectures:

\begin{conj}[Absolute Suspension Conjecture]\label{ASC}
Every $n$-dimensional absolute suspension is homeomorphic to
the $n$-sphere.
\end{conj}

\begin{conj}[Absolute Suspension Conjecture]\label{ACC}
Every $n$-dimensional absolute cone is homeomorphic to the $n$-cell.
\end{conj}

The fact that small metric balls in Busemann $G$-spaces are absolute
cones follows from the local cone structure and homogeneity.
Therefore the truth of the Absolute Cone Conjecture~\ref{ACC}
would imply the truth
of the Busemann Conjecture~\ref{BC}.

In 1974 Szyma\'{n}ski \cite{Szymanski} proved the Absolute
Suspension Conjecture~\ref{ASC}
for dimensions $n \leq 3$.  In 1978, Mitchell
\cite{Mitchell1} gave an alternate proof to Szyma\'{n}ski result and
showed that every $n$-dimensional absolute suspension is an ENR
homology $n$-manifold homotopy equivalent to the $n$-sphere. 
In 2005 Bellamy and Lysko \cite{Bellamy-Lysko}
proved the generalized Sch\"{o}nflies theorem for absolute suspensions.
and
Nadler \cite{Nadler} gave a proof of the Absolute Cone
Conjecture~\ref{ACC} in dimensions $n=1,2$. However, in 2007 Guilbault
\cite{Guilbault} completely clarified the status of the Absolute Cone
Conjecture~\ref{ACC}:

\begin{thm}[Guilbault]
The Absolute Cone
Conjecture~\ref{ACC} is true for $n \leq 4$ and false for $n \geq 5$.
\end{thm}

  Guilbaut proved this result in dimensions $n \geq 5 $ by
constructing counter-examples. For the special case $n=4$, Guilbaut
shows the Absolute Cone
Conjecture~\ref{ACC}
is true modulo the
$3$-dimensional Poincar\'{e} Conjecture, which indeed follows by
Perelman's proof of the Poincar\'{e} Conjecture in dimension $n=3$
\cite{Morgan-Tian}. Although the solution of the Absolute Cone
Conjecture~\ref{ACC}
leaves the status of the Busemann Conjecture~\ref{BC}
unresolved, it does cast some
suspicion the validity of \ref{BC}.

\subsection{The Moore Conjecture}

A space $X$ is said to be a {\it codimension one manifold  factor}
if $X \times \mathbb{R}$ is a topological manifold.  In 1955 Bing
constructed his infamous Dogbone space \cite{Bing1}.
Bing's Dogbone space is the
image of a cell-like map $\pi: \mathbb{R}^3 \to X$.  Bing \cite{Bing2}
showed
that the the Dogbone space $X$ is not homeomorphic to
$\mathbb{R}^3$, however $X \times \mathbb{R}^1$ is homeomorphic to
$\mathbb{R}^4$. This result led to the Moore Conjecture:

\begin{conj}[Moore Conjecture]\label{MC}
Every resolvable generalized manifold is a codimension one manifold
factor.
\end{conj}

The Moore Conjecture~\ref{MC} is also related to the Busemann Conjecture~\ref{BC}.
Every Busemann $G$-space is a manifold if and only if small metric
spheres are codimension one manifold factors. Equivalently, in
dimensions $n \geq 5$, every Busemann $G$-space $X$ is a manifold if
and only if small metric spheres $\Sigma \subset X$ are resolvable
and have the property that $\Sigma \times \mathbb{R}$ has the
disjoint disks property.

Although it is unknown whether small metric spheres are resolvable,
Thurston showed that they are generalized $(n-1)$-manifolds. Also,
according to the properties of the Quinn index number which measures
the obstruction of a space being resolvable, the resolvability of
$\Sigma$ is equivalent to the resolvability of $X$ (see
\cite{Quinn2}). Moreover, Mitchell \cite{Mitchell1} proved in 1978
that any $n$-dimensional absolute suspension $X$ is a regular
generalized $n$-manifold homotopy equivalent to $S^n$;  all its
links are generalized $(n-1)$-manifolds homotopy equivalent to
$S^{n-1}$.

He furthermore showed that an $n$-dimensional absolute
cone $X$ is a regular generalized $n$-manifold proper homotopy
equivalent to $\mathbb{R}^n$;  all its links are generalized
$(n-1)$-manifolds homotopy equivalent to $S^{n-1}$.  Note that if in
Mitchell's theorem "homotopy equivalent" could be replaced with
"fine homotopy equivalent", Mitchell's theorem would imply
resolvability \cite{CF} (see also \cite{Lacher}\cite{Mitchell2}\cite{Mitchell-Repovs}).

Although it is also still unknown whether small metric spheres $X$
in Busemann $G$-spaces satisfy the disjoint disks properties, there
have been several results determining useful general position
properties of an ANR $X$ that characterize $X \times \mathbb{R}$ as
having the the disjoint disks property.  In particular, these
properties include:

    (i) {\it The disjoint arc-disk property}
    (Daverman \cite{Daverman}).
    A space $X$ has the disjoint arc-disk property
    if any pair of maps $\alpha: I \to X$ and $f: D^2 \to X$ can be approximated by
    paths with disjoint images (i.e., $X$ has the disjoint $(1,2)$-cells property).  If $X$ has the disjoint arc-disk
    property, then $X \times \mathbb{R}$ has the disjoint disks
    property.

    (ii) {\it The disjoint homotopies property}
    (Edwards \cite{Edwards}, Halverson
    \cite{Halverson1}).
    A space $X$ has the disjoint homotopies
    property if every pair of path homotopies $f,g: D
    \times I \to X$, where $D=I=[0,1]$, can be approximated by homotopies $f',g': D
    \times I \to X$ so that $f_t(D) \cap g_t(D) = \emptyset$ for all
    $t \in I$.  If $X$ has the disjoint homotopies
    property, then $X \times \mathbb{R}$ has the disjoint disks
    property.

    (iii) {\it The plentiful $2$-manifolds property}
    (Halverson
    \cite{Halverson1}).
    An ANR $X$ has the plentiful $2$-manifolds
    property if every path $\alpha: I \to X$ can be approximated by
    a path $\alpha': I \to N\subset X$ where $N$ is a $2$-manifold
    embedded in $X$.  If an ANR $X$ of dimension $n \geq 4$ has the
    plentiful $2$-manifolds property, the $X$ has the disjoint
    homotopies property.

    (iv) {\it The method of $\delta$-fractured maps}
    (Halverson
    \cite{Halverson2}).
    A map $f:D \times I \to X$ is said to be $\delta$-fractured over a
    map $g:D \times I \to X$, where $D = I = [0,1]$, if there are disjoint balls $B_1, B_2,
    \ldots, B_m$ in $D \times I$ such that
    \begin{enumerate}
        \item $diam(B_i) < \delta$;
        \item $f^{-1}(im(g)) \subset \bigcup_{i=1}^m int(B_i)$; and
        \item $diam(g^{-1}(f(B_i))) < \delta$
    \end{enumerate}
    If $X$ is a space such that for any two path homotopies $f,g: D
    \times I \to X$, where $g$ is a constant path homotopy, and
    $\delta >0$, there are approximations $f', g': D \times I \to X$
    such that $f'$ is $\delta$-fractured over $g'$, then $X$ has the
    disjoint homotopies property.

    (v) {\it The $0$-stitched disks property}
    (Halverson
    \cite{Halverson3}). A space $X$ has the \emph{$0$-stitched disks property} if any two
    maps $f,g:D^2 \to X$ can be approximated  by maps $f',g':D^2 \to X$
    such that there are infinite $1$-skeleta
    $(K^\infty)^{(1)}$ and $(L^\infty)^{(1)}$ of $D^2$ and
    $0$-dimensional $F_\sigma$ sets $A \subset \text{int}(D^2) -
    (K^\infty)^{(1)}$ and $B \subset  \text{int}(D^2) -
    (L^\infty)^{(1)}$ such that $f'|_{(K_1^\infty)^{(1)}} \cup
    g'|_{(K_2^\infty)^{(1)}}$ is $1-1$ and $f'(D^2-A) \cap g'(D^2-B) =
    \emptyset$. If $X$ has the $0$-stitched disks property, then $X$ has
    DHP.

    (vi) {\it The disjoint concordances property}
    (Daverman and Halverson
    \cite{DaH}).  A \emph{path concordance} in a space $X$ is a
    map $F:D \times I \to X \times I$, where $D=I=[0,1]$, such that $F(D
    \times e) \subset X \times e, e \in \{0,1\}.$  A metric space
    ($X,\rho$) satisfies the \emph{Disjoint Path Concordances Property
    (DCP)} if, for any two path homotopies $F_i:D \times I \to X$
    ($i=1,2$) and any $\epsilon > 0$, there exist path concordances
    $F'_i: D \times I \to X \times I$ such that
    \begin{center}
    $F'_1(D \times I) \cap F'_2(D \times I) = \emptyset$
    \end{center}
    and $\rho (F_i, \text{proj}_X F'_i) < \epsilon$.  An ANR $X$ has
    the disjoint concordances property if and only if $X \times
    \mathbb{R}$ has the disjoint disks property.

Due to homogeneity, if a Busemann $G$-space $X$ has a single metric
sphere satisfying any one of these properties, then $X$ has the
disjoint disks property.

\section{Summary and Questions}

In summary the following relationships hold between the conjectures
and problems discussed in this survey.
\begin{itemize}
\item {Bing-Borsuk Conjecture~\ref{BBC}
      $\Rightarrow$
      Busemann Conjecture~\ref{BC}}
\item {de Groot Conjecture~\ref{ACC}
      $\Rightarrow$
      Busemann Conjecture~\ref{BC}}
\item {Moore Conjecture~\ref{MC} \&
      Busemann $G$-Spaces Resolution Conjecture~\ref{BSRC}}
      $\Rightarrow$
      {Busemann Conjecture~\ref{BC}}
\item {Bryant-Ferry-Mio-Weinberger Conjecture~\ref{BFMW1}}
      $\Rightarrow$
      {failure of the Bing-Borsuk Conjecture~\ref{BBC}}
\end{itemize}

So far, the validity of only one of these implications has been
determined.  Recall that the de Groot Conjecture~\ref{ACC} was
shown to be false for all $n\ge 5$ (see \cite{Guilbault}).  Note
that the failure of the Busemann Conjecture~\ref{BC} would settle
the Bing-Borsuk Conjecture~\ref{BBC} in the negative. If the
Busemann $G$-Spaces Resolution Conjecture~\ref{BSRC} were proved and the
Busemann Conjecture~\ref{BC} proved to be false, then the Moore
Conjecture~\ref{MC} would be
settled in the negative.

Below is a summary of relevant questions that remain unsolved:

\begin{enumerate}
\item Are all Busemann $G$-spaces resolvable?
\item Do all Busemann $G$-spaces $X$ of dimension $n \geq 5$
have the disjoint disks property? (or equivalently, does $X$ contain
some metric sphere $\Sigma$ that has a general position property
that implies $X \times \mathbb{R}$ has DDP)?
\item Are all finite-dimensional Busemann $G$-spaces manifolds?
\item Are all absolute cones resolvable?
\item Are all finite-dimensional homogeneous connected compact metric spaces resolvable?
\item Are all resolvable generalized manifolds codimension one
manifold factors?
\item Are all generalized manifolds with the disjoint disks property
homogeneous?
\end{enumerate}

\section{Epilogue: Homogeneity and group actions}

The Bing-Borsuk Conjecture~\ref{BBC} belongs to a wide group of difficult open
problems related to homogeneity and group actions. The nearest one
is an old problem: Is the Hilbert cube is the only homogeneous
compact AR?

All problems of this sort can be seen in the following framework:
Given a topological group $G$ and a closed subgroup $H$, describe
the topological structure of the coset space $G/H$ assuming
that it has some extra properties (local contractibility,
finite-dimensionality, local compactness, etc.) Model results here
concern the structure of topological groups \cite{DoT}\cite{MZ}:

\begin{thm}[Montgomery-Zippin]
Each locally compact locally contractible topological group is a Lie group
and hence a manifold.
\end{thm}

A {\it Polish group}  is a topological group which is also a
Polish space \cite{BK}:

\begin{thm}[Dobrowolski-Torunczyk]
Each Polish ANR-group is a Hilbert manifold (finite or infinite-dimensional).
\end{thm}

The last theorem suggests the following problem which was partially solved by Banakh and Zarichnyi \cite{BZ}:

\begin{prob}
Is each complete metric ANR-group a Hilbert manifold?
\end{prob}

Now let us turn to homogeneous spaces.

\begin{prob}[Banakh]\label{main}
Let $G$ be a Polish group and $H$ a closed subgroup such that $G/H$
is an ANR. Is then $G/H$ a manifold modeled on: (i) the Euclidean
$n$-space $\mathbb{R}^n$; (ii) the Hilbert cube $Q$; or
(iii) the Hilbert sequence  space $l_2$? \\
What if $G$ is an ANR-group?
What if the quotient map $G\to
G/H$ is a locally trivial bundle?
\end{prob}

Problem~\ref{main}(i) is exactly the Bing-Borsuk Problem
while the second part is related to the  question mentioned above
on homogeneous
compact AR's. For connected locally compact topological groups,
Problem~\ref{main} was answered by Szenthe \cite{Szenthe1}\cite{Szenthe}.

\begin{thm}[Szenthe]\label{szenthe}
Let $G$ be a locally compact topological group such that
the quotient group $G/G_0$ of $G$ by its  (connected) identity component $G_0$ is compact
(the quotient group $G/G_0$ is called the group of components and denoted $\pi_0(G)$.)
Then for any closed subgroup $H\subset G$, the coset space $G/H$ is a
disjoint union of
topological manifolds if and only if it is locally contractible.
\end{thm}

This theorem has some interesting consequences for homogeneous metric spaces.
We define a metric space $X$ to be {\em metrically homogeneous} if for any two points $x,y\in X$ there is an isometry $f:X\to X$ such that $f(x)=y$. Note that every  $C^{\infty}$-smooth Riemannian manifold $M$ is Busemann G-space, hence $M$ is homogeneous. However, $M$ is in general not metrically homogeneous. We also observe that every transitive group $G$
of isometies of every locally compact metric space
$M$ admits a natural metric with respect to which
$G$ is locally compact and acts continuously on $M$.
It follows from Theorem~\ref{szenthe} that the
"isometric version" of the Bing-Borsuk Conjecture~\ref{BBC} is true:

\begin{cor}\label{c1} A metrically homogeneous compact metric space $X$ is a topological
$n$-manifold if and only if $X$ is locally compact and locally contractible.
\end{cor}

Results from \cite{Szenthe1}\cite{Szenthe} also imply that any locally compact
connected (possibly metrizable) locally contractible
topological space $M$ whith a locally compact
transitive continuous group $G$ of homeomorphisms,
necessarily admits a structure of a
$C^{\infty}$-manifold
and a compatible Riemannian metric tensor $g$ such
that $G$ acts by isometries on $(M,g).$

This implies two facts:
(i)
such a manifold $M$ is necessarily smoothable. As a corollary,
no nonsmoothable compact $4$-manifold (recall that most compact
4-manifolds are nonsmoothable \cite{Freedman-Quinn}\cite{Scorpan}) admits any locally compact
continuous transitive group of homeomorphisms; and
(ii)
there are many smooth manifolds which admit no metrically
homogeneous Riemannian metric (for example, the $2$-sphere
with two handles). Therefore no such manifold admits any locally
compact transitive continuous group of homeomorphisms.

A topological  space $X$ is said to be
{\em continuously homogeneous}
if for every $x,y\in X$ there is a homeomorphism
 $h_{x,y}:X\to X$
such that
$h_{x,y}(x)=y$ and
$h_{x,y}$
continuously
depends on the points $x,y$ in the sense that the map
$H:X^3\to X$,
$H:(x,y,z)\mapsto h_{x,y}(z)$, is continuous.

It is easy to see that
each topological group (endowed with a left-invariant metric) is
continuously homogeneous (as a metric space). Continuously
homogeneous spaces were introduced and studied by Banakh et al.
\cite{BKT2}\cite{BKT}. It can be shown that a topological space $X$
is continuously homogeneous if and only if it is rectifiable in the
sense of Gulko \cite{Gul}).

\begin{prob}[Banakh]\label{pr3}
Let $X$ be a continuously homogeneous Polish ANR-space. Is $X$ a Hilbert manifold?
\end{prob}

It light of this problem one should mention that continuously homogeneous spaces
cannot be Hilbert cube manifolds \cite{BKT} (which distinguishes Problem~\ref{pr3}
from a more general Problem~\ref{main}).

\section{Acknowledgements}
The authors acknowledge several comments and suggestions from T. O. Banakh and V. N. Berestovskii.
The authors were supported in part by the Slovenian Research Agency grants
BI-US/07-08/029,
P1-0292-0101, and
J1-9643-0101.
This paper was presented at the 4th Croatian Mathematical Congress (Osijek, June 17-20, 2008).
We thank the organizers for the invitation and hospitality.

\end{document}